\documentclass{amsart}

\usepackage{amssymb,amscd,verbatim,url}
\usepackage{enumerate}

\long\def\citat#1{{\narrower\narrower\small  #1 \par \noindent}}

\usepackage{local_KK80}   

\author
{Bernhelm Boo{\ss}-Bavnbek}
\address{NSM/IMFUFA, Roskilde University, P.O. Box 260,
DK-4000 Roskilde, Denmark}
\email{booss@ruc.dk}

\begin{document}


\title
{The Mathematization of the Individual Sciences - Revisited}




\begin{abstract}
We recall major findings of a systematic investigation of the mathematization
of the individual sciences, conducted by the author
in Bielefeld some 35 years ago under the direction of Klaus Krickeberg, and
confront them with recent developments in physics, medicine, economics, and spectral geometry.
\end{abstract}

\maketitle

\par
\centerline{\bf Dedicated to Klaus Krickeberg on his 80th birthday}\par


\section*{Introduction - recollections}

In the years 1972-75, a small working group on {\em The
mathematization of the individual sciences} was organized by Klaus
Krickeberg and myself at the then young West-German university of
Bielefeld. The group included representatives of all the disciplines
represented at the university. Unfortunately, the engineering and
the medical sciences were not represented, so that the important
questions of mathematization of industrial production and health
care were not addressed.

To my knowledge it was the first systematic investigation of the
mathematization as seen from the individual sciences. The term {\em
mathematization} was directed both to the shining modelling aspects
(= selecting, finding, inventing the {\em right} specification,
equations and set-up) and the rather profane tasks of executing the
necessary calculations (= parameter estimations, analytic solutions,
numerical simulations, approximations,
stability/deformation/robustness arguments, and geometric
interpretations). The investigation remained a singularity in its
multi-disciplinary, almost all-embracing ambition. It was based on
extensive literature study and a long series of hearings and
interviews. It was, and perhaps continues to be an offence against
two traditional imperatives of scientific research: {\em Don't
mix}\/! and {\em Originality first}\/! We were continuously
comparing - and confronting - different sciences with each other and with
their mathematization experiences, looking forth and back, mixing
continuously mathematical, methodological, epistemological,
educational, and ethical questions. More offending, we were not
interested in the most recent top-notch results, but in the {\em
relevance} of the applied mathematical methods for advances in the
individual sciences, and in the question of what the cases were {\em
representative} for.

On the initiative of
the wizard Alexander Ostrowski (1893-1986),
Krickeberg and I collected some of our findings in a volume
and published it with Birkh{\"a}user in Basel in 1976, \cite{BooKri:MEW}, see also \cite{Boo:NTM}
for a summary report.

\subsection*{Acknowledgment} For the following recollections and generalizations, I
must take the responsibility alone, assuming that most of the views
were or will be shared by Klaus Krickeberg. The findings are also
based on later modelling experiences and continuing discussions with
students and colleagues there. In particular, I am indebted to Viggo
Andreasen, Peder Voetmann Christiansen, Niels Lang\-ager Ellegaard,
Jens H{\o}yrup, Jens H{\o}jgaard Jensen, Bent C. J{\o}rgensen,
Jesper Larsen, Anders Madsen, Mogens Niss, Johnny Ottesen, and Stig
Andur Pedersen (all Roskilde), Philip J. Davis (Providence),
Giampiero Esposito (Napoli), Roman Galar (Wroclaw), Martin Koch
(Copenhagen), Matthias Lesch (Bonn), Glen Pate (Hamburg) and Erik
Renstr{\"o}m (Lund).

\section{Continuing challenges}\label{s:cc}
The use of mathematical arguments, first in pre-scientific investigations,
then in other sciences,
foremost in medicine and astronomy and in their shared border region astrology, has been
traced way back in history by many authors from various perspectives,
Bernal \cite{Ber:SIH}, H{\o}yrup \cite{Hoy:FHS}, and Kline \cite{Kli:MTF}.

Globally speaking, they all agree on three mathematization tendencies:
\begin{enumerate}
\item The progress in the individual sciences makes work on ever more complicated problems possible and necessary.
\item This accumulation of problems and data demands conscious, planned, and economic
procedures in the individual sciences, i.e., an increased emphasis on questions of
methodology.
\item Finally, this increased emphasis on questions of methodology is as a rule associated
with the tendency of mathematization.
\end{enumerate}

All of this applies generally. In detail, we find many various
pictures. In his {\em Groundplat of Sciences and Artes,
Mathematicall} of 1570 \cite{Dee:GSA}, the English alchemist,
astrologer, and mathematician John Dee, the first man to defend the
Copernican theory in Britain and a consultant on navigation, pointed
out, in best Aristotelian tradition, that it is necessary in the
evaluation of mathematization to pay strict attention to the
specific characteristics of the application area in question. He
postulated a dichotomy between the {\em Principall} side, pure
mathematics, and the {\em Deriuative} side, i.e., applied
mathematics and mathematization. He then classified the applications
of pure mathematics according to objects treated:
\begin{itemize}
\item {\em Ascending Application in thinges Supernaturall, eternall and Diuine},
\item {\em In thinges Mathematicall: Without farther Applications},
\end{itemize}
and finally, on the lowest and most vulgar plane in the Aristotelian scheme,
\begin{itemize}
\item {\em Descending Application in thinges Naturall: both Substantiall \& Accidentall, Visible \& Inuisible \& c.}.
\end{itemize}
Now that history has excluded matters divine from mathematics, we can with some justification ask whether
later generations may regard with equal amusement and astonishment the fact that in our time there are
a large number of professional mathematicians, who are completely satisfied with spending their entire lives
working in the second, inner mathematical level and who persistently refuse to descend to vulgar applications.

The panorama of the individual sciences and the role that mathematics had to play in them was perfectly
clear for John Dee. In our time the matter is somewhat more complex. In this review, I cannot
point out a geodetically perfect picture of today's landscape of mathematization. I must treat the matter
rather summarily. A summary treatment may have the advantage that in comparison among the mathematization progress
in three sciences below (physics, medicine and economics), common problems on one hand and
special features of mathematization on the other hand can be seen more clearly.

In the following, I shall ``ascend" from the study of dead nature in
{\em physics}, the field which has the highest degree of
mathematization on any chosen scale, both quantitatively and
qualitatively, over the investigation of living matter in {\em
medicine}, the field where one might expect the greatest
mathematization advances in our century, to the treatment of
financial issues and decision making for commerce and production in
{\em economics}, a field of questionable scientific state, that,
beyond well founded actuary estimations, lacks unambiguous results
and convincing clear perspectives regarding mathematization.

\section{Physics}\label{s:physics}
The intimate connection between mathematics and physics makes it
difficult to determine the theoretical relevance of mathematics and
obscures the boundary between genuinely physical thought and
observation on one side and the characteristically mathematical
contribution on the other side. Recall Hilbert's perception of
probability theory as a chapter of physics in his famous 6th Problem
\cite{Hil:MP}:

\citat{6. Mathematical Treatment of the Axioms of Physics. The
investigations on the foundations of geometry suggest the problem:
{\em To treat in the same manner, by means of axioms, those physical
sciences in which already today mathematics plays an important part;
in the first rank are the theory of probabilities and mechanics}.}

\noi To say it mildly, as Gnedenko did in his comments to the Russian edition of 1969: {\em Today
this viewpoint} (to consider probability theory as a chapter of physics) {\em
is no longer so common as it was around the turn of the century, since the independent mathematical content of
the theory of probabilities has sufficiently clearly showed since then...}
With hindsight and in view of the still challenging foundational problems of quantum mechanics,
however, we may accept that
parts of mathematics and physics can be interlaced in a non-separable way.

Another famous example of that inextricable interlacement is
provided by the Peierls-Frisch memorandum of 1940 to the British
Government: suggested by the codiscoverer of fission Otto Frisch,
the physicist Rudolph Peierls, like Frisch a refugee in Britain,
made the decisive feasibility calculation that not tons (as -
happily - erroneously estimated by Heisenberg in the service of the
Nazis) but only about 1 kg (later corrected to 6 kg) of the pure
fissile isotope $U_{235}$ would be needed to make the atomic bomb.
Was it mathematics or physics? It may be worth mentioning that
Peierls was a full professor at the University of Birmingham since
1937 and became joint head of {\em mathematics} there,
\cite{Edw:REP}. Theoretical physics in Britain is often in
mathematics. As a matter of fact, physics in our sense did not exist
as a single science before the nineteenth century. There were
well-defined {\em experimental physics} comprising heat, magnetism,
electricity and colour, leaving mechanics in mathematics, see
\cite[p. 493]{Hoy:FHS}.

In spite of that intermingling, physics can provide a filter for our
review, a ready system of categories to distinguish different use of
mathematics in different modelling situations. Perhaps, the
situation can be best compared with the role of physics in general
education: after all, physics appears as the model of
mathematization: there is no physics without mathematics - and, as a
matter of fact, learning of mathematics is most easy in a physics
context: calculation by letters; the various concepts of a function
(table, graph, operation) and its derivatives and anti-derivatives;
differential equations; the concept of observational errors and the
corresponding estimations and tests of hypotheses; Brownian
movements; all these concepts can be explained context-free or in
other contexts (where some of the concepts actually originated), but
they become clearest in the ideally simple applications of physics,
which are sufficiently complicated to see the superiority of the
mathematization as compared to feelings, qualitative arguments,
discussions, convictions, imagination - but simple enough to get
through.

\subsection{Variety of modelling purposes}\label{ss:variety}
It may be helpful to distinguish the following modelling purposes:

\subsubsection{Production of data, model based measurements}
Clearly, the public associates the value of mathematical modelling
foremost to its {\em predictive} power, e.g., in numerical weather
prediction, and its {\em prescriptive} power, e.g., in the design of
the internal ballistics of the hydrogen bomb; more flattering to
mathematicians, the {\em explanatory} power of mathematization and
its contribution to {\em theory development} yield the highest
reputation within the field. However, to the progress of physics,
the {\em descriptive} role, i.e., supporting model-based
measurements in the laboratory, is -- as hitherto -- the most
decisive contribution of mathematics. Visco-elastic constants and
phase transition processes of glasses and other soft
materials can not be measured directly. For high precision in the
critical region, one measures electric currents through a ``dancing"
piezoelectric disc with fixed potential and varying frequency. In
this case, solving mathematical equations from the fields of
electro-dynamics and thermo-elasticity becomes mandatory for the
design of the experiments and the interpretation of the data. In
popular terms, one may speak of a {\em mathematical microscope}, in
technical terms of a {\em transducer} that becomes useful as soon as
we understand the underlying mathematical equations.

\subsubsection{Simulation}\label{ss:simulation}
Once a model is found and verified and the system's parameters are
estimated for one domain, one has the hope of doing computer
calculations to predict what new experiments in new domains (new
materials, new temperatures etc) should see. Rightly, one has given
that type of calculations a special name of honour, {\em computer
simulations}: as a rule, it requires to run the process on a
computer or a network of computers under quite sophisticated
conditions: typically, the problem is to bring the small distances
and time intervals of well-understood molecular dynamics up to
reasonable macroscopic scales, either by aggregation or by Monte
Carlo methods -- as demonstrated by Buffon's needle casting for the
numerical approximation of $\pi$.

One should be aware that the word ``simulation" has, for good and
bad, a connation derived from NASA's space simulators and Nintendo's
war games and juke boxes. Animations and other advanced computer
simulations can display an impressive beauty and convincing power.
That beauty, however, is often their dark side: simulations can show
a deceptive similarity with true observations, so in computational
fluid dynamics when the numerical solution of the Bernoulli
equations, i.e., the linearization of the Navier-Stokes equations
for laminar flow displays eddies characteristic for the non-linear
flow. The eddies do not originate from real energy loss due to
friction and viscosity but from hardly controllable hardware and
software properties, the chopping of digits, thus providing a
{\em  magic realism}, as coined by Abbott and Larsen
\cite{AbbLar:MCD}. In numerical simulation, like in mathematical
statistics, results which fit our expectations too nicely, must
awake our vigilance instead of being taken as confirmation.


\subsubsection{Prediction}\label{ss:prediction}
As shown in the preceding subsection, there is no sharp boundary
between description and prediction. However, the quality criteria
for predictions are quite simple: do things develop and show up as
predicted? So, for high precision astrology and longitudinal
determination in deep-sea shipping the astronomical tables of
planetary movement, based on the outdated and falsified Ptolemaic
system (the {\em Resolved Alfonsine Tables}) and only modestly
corrected in the {\em Prutenic Tables} of 1551 were, until the
middle of the 17th century, rightly considered as more reliable than
Kepler's heliocentric {\em Rudolphine Tables}, as long as they were
more precise - no matter on what basis, see, Steele \cite[p.
128]{Ste:OPE}.

Almost unnoticed, we have had a similar revolution in weather
prediction in recent years: the (i) analogy ({\it synoptic}) methods of identifying a
similarly looking weather situation in the weather card archives to
base the extrapolation on it were replaced by almost pure (ii)
numerical methods to derive the prediction solely from the
thermodynamic and hydrodynamic basic equations and conservation
laws, applied to initial conditions extracted from the observation
grid. ``Almost" because the uncertainty of the interpolation of the
grid and the high sensitivity of the evolution equations to initial
conditions obliges to repeated runs with small perturbations and
human inspection and selection of the most ``probable" outcome like
in (i). That yields sharp estimates about the certainty of the
prediction for a range of up to 10 days. In nine of ten cases, the
predictions are surprisingly reliable and would have been impossible
to obtain by traditional methods. However, a 10{\%} failure rate
would be considered unacceptable in industrial quality
control.

In elementary particle physics, the coincidence of predictions with
measurements is impressive, but also disturbing. I quote from Smolin
\cite[pp. 12-13]{Smo:TWP}:

\citat{Twelve particles and four forces are all we need to
explain everything in the known world. We also understand very well
the basic physics of these particles and forces. This understanding
is expressed in terms of a theory that accounts for all these
particles and all of the forces except for gravity. It's called {\em
the standard model of elementary-particle physics} - or the standard
model for short. \dots Anything we want to compute in this theory we
can, and it results in a finite number. In the more than thirty
years since it was formulated, many predictions made by this theory
have been checked experimentally. In each and every case, the theory
has been confirmed.

The standard model was formulated in the early 1970s.
Except for the discovery that neutrinos have mass, it has not required adjustment since.
So why wasn't physics over by 1975? What remained to be done?

For all its usefulness, the standard model has a big problem: It has a long list of adjustable constants. \dots}

\noi We feel pushed back to the pre-Keplerian, pre-Galilean and pre-Newtonian cosmology built on
ad-hoc assumptions, displaying clever and deceptive mathematics-based similarity between observations
and calculations -- and ready to fall at any time because the basic assumptions are not explained.

Perhaps the word {\em deceptive} is inappropriate when speaking of
description, simulation and prediction: for these tasks, {\em
similarity} can rightly be considered as the highest value
obtainable, as long one stays in a basically familiar context. From
a semiotic angle, the very similarity must have a meaning and is
indicating something; from a practical angle, questions regarding
the epistemological status can often be discarded as metaphysical
exaggerations: who cares about the theoretical or ad-hoc basis of a
time schedule in public transportation -- as long as the train leaves
on time!

\subsubsection{Control}\label{ss:control}
The prescriptive power of mathematization deserves a more critical examination.
In physics and engineering we may distinguish between the (a) feasibility, the
(b) efficiency, and the (c) safety of a design. A {\em design} can be an object
like an airplane or a circuit diagram for a chip, an instrument like a digital
thermometer, TV set, GPS receiver or pacemaker, or a regulated process like a
feed-back regulation of the heat in a building, the control of a power station or the
precise steering of a radiation canon in breast cancer therapy.
Mathematics has its firm footing for testing (a) in thought experiments, estimations of process parameters,
simulations and solving equations. For testing (b), a huge inventory is available of
mathematical quality control and optimisation procedures by variation of key parameters.

It seems to me, however, that (c), i.e., safety questions provide
the greatest mathematical challenges. They appear differently in (i)
experience-based, (ii) science-based and (iii) science-integrated
design. In (i), mathematics enters mostly in the certification of
the correctness of the design copy and the quality test of the
performance. In (ii), well-established models and procedures have to
be modified and re-calculated for a specific application.
Experienced physicists and engineers, however, seldom trust their
calculations and adaptations. Too many parameters may be unknown and
pop-up later: Therefore, in traditional railroad construction, a
small bridge was easily calculated and built, but then
photogrammetrically checked when removing the support constructions.
A clash of more than $\delta_{\operatorname{crit}}$ required
re-building. Similarly, even the most carefully calculated chemical
reactors and other containers under pressure and heat have their
prescribed ``Soll-Bruchstelle" (supposed line of fracture) in case
that something is going wrong.

The transition from (ii) to (iii) is the most challenging: very
seldom one introduces a radically new design in the physics
laboratory or engineering endeavour. But there are systems where all
components and functions can be tested separately but the system as
a whole can only be tested {\em in situ}: a new design of a Diesel
ship engine; a car, air plane or space craft; a new concept in
cryptography. In all these cases, one is tempted to look and even to
advocate for mathematical proofs of the safe function according to
specification. Unfortunately, in most cases these ``proofs" belong
rather to the field of fiction than to rigorous mathematics. For an
interesting discussion on ``proofs" in cryptography (a little remote
from physics) see the debate between Koblitz and opponents in
\cite{Kob:xxx} and follow-ups in the {\em Notices of the American
Mathematical Society}.

An additional disturbing aspect of science-integrated technology
development is the danger of a loss of transparency. Personally, I
must admit, I'm grateful for most black-box systems. I have no
reason to complain when something in my computer is hidden for my
eyes, as long everything functions as it shall or can easily be
re-tuned. However, for the neighbourhood of a chemical plant (and
the reputation of the company) it may be better not to automatize
everything but to keep some aspects of the control non-mathematized
and in the hands of the service crew to avoid de-qualification and
to keep the crew able to handle non-predictable situations.

A last important aspect of the prescriptive power of the
mathematization is its formatting power for thought structure and
social behaviour. It seems that there is not so much to do about it
besides being aware of the effects.

\subsubsection{Explain phenomena}
The noblest role of mathematical concepts in physics is to
explain phenomena. Einstein did it when {\em reducing} the heat
conduction to molecular diffusion, starting from the formal analogy
of Fick's Law with the cross section of Brownian motion. He did it
also when {\em generalizing} the Newtonian mechanics into the
special relativity of constant light velocity and again when {\em
unifying} forces and curvature in general relativity.

Roughly speaking, mathematical models can serve physics by reducing
new phenomena to established principles; as heuristic devices for
suitable generalizations and extensions; and as ``a conceptual
scheme in which the insights \dots fit together" (C. Rovelli).
Further below we shall return to the last aspect -- the unification
hope.

Physics history has not always attributed the best credentials to
explaining phenomena by abstract constructions. It has discarded the
concept of a ghost for perfect explanation of midnight noise in old
castles; the concept of ether for explaining the finite light
velocity; the phlogiston for burning and reduction processes, the
Ptolemaic epicycles for planetary motion. It will be interesting to
see in the years to come whether the mathematically advanced
String Theory or the recent Connes-Marcolli reformulation of the
Standard Model in terms of spectral triples will undergo the same
fate.

\subsubsection{Theory development}
Finally, what has been the role of mathematical concepts and
mathematical beauty for the very theory development in physics? One
example is Johann Bernoulli's purely aesthetic confirmation of
Galilean fall law $s\,=\, g/2 \ t^2$ among a couple of candidates as
being the only one providing the same equation (shape) for his
brachistochrone and Huygens' tautochrone, \cite[p. 395]{Ber:xxx}:

\citat{Before I end
I must voice once more the admiration that I feel for the unexpected identity of Huygens'
tautochrone and my brachystochrone. I consider it especially remarkable that this coincidence
can take place only under the hypothesis of Galilei, so that we even obtain from this
a proof of its correctness. Nature always tends to act in the simplest way, and so it here lets
one curve serve two different functions, while under any other hypothesis we should need
two curves.}

Another, more prominent example is the lasting triumph of Maxwell's equations: a world of radically
new applications were streaming out of the beauty and simplicity of the equations of electro-magnetic waves!

However, not every mathematical, theoretical and empirical
accumulation leads to theory development. Immediately after
discovering the high-speed rotation of the Earth around its own
axis, a spindle shape of the Earth was suggested and an
infinitesimal tapering towards the North pole confirmed in geodetic
measurements around Paris. Afterwards, careful control measurements
of the gravitation at the North Cap and at the Equator suggested the
opposite, namely an ellipsoid shape with flattened poles. Ingenious
mathematical mechanics provided a rigorous reason for that. Gauss
and his collaborator Listing, however, found something different in
their control. They called the shape {\em gleichsam
wellenf{\"o}rmig} and dropped the idea of a theoretically
satisfactory description. Since then we speak of a {\em Geoid}.

\subsection{``The trouble with physics"} That is the title of an interesting and well-informed polemic
by Lee Smolin against String Theory and present main stream physics at large. He notices a {\em stagnation}
in physics, {\em so much promise, so little fulfillment} \cite[p. 313]{Smo:TWP}, a predominance of
anti-foundational spirit and contempt for visions, partly related to the mathematization paradigm of the 1970s,
according to Smolin:  {\em Shut up and calculate}.

It seems to me that Smolin, basically, may be right. B{\o}rge
Jessen, the Copenhagen mathematician and close collaborator of
Harald Bohr once suggested to distinguish in sciences and
mathematics between periods of {\em expansion} and periods of {\em
consolidation}. Clearly physics had a consolidation period in the
first half of the 20th century with relativity and quantum mechanics.
The same may be true for biology with the momentous triumph of the
DNA disclosure around 1950, while, to me, the mathematics of that
period is characterized by an almost chaotic expansion in thousands of
directions. Following that way of looking, mathematics of the second
half of the 20th century is characterized by an enormous consolidation,
combining so disparate fields like partial differential equations
and topology in index theory, integral geometry and probability in
point processes, number theory, statistical mechanics and
cryptography, etc. etc. A true period of consolidation for
mathematics, while - at least from the outside - one can have the
impression that physics and biology of the second half of the 20th
century were characterized merely by expansion, new measurements,
new effects - and almost total absence of consolidation or, at least
failures and vanity of all trials in that direction.

Indeed, there have been impressive successes in recent physics,
in spite of the absence of substantial theoretical progress in physics: perhaps the most spectacular
and for applications most important discovery has been the High Temperature Superconducting (HTS) property
of various ceramic materials by Bednorz and M{\"u}ller - seemingly without mathematical or theoretical efforts
but only by systematic combinatorial variation of experiments - in the tradition of the old alchemists, \cite{BedMul:HTS}.

The remarkable advances in fluid dynamics, weather prediction,
oceanography, climatic modelling are mainly related to new
observations and advances in computer power while the equations have
been studied long before.

Nevertheless, I noticed a turn to theory among young experimental physicists in recent years, partly related to
investigating the {\em energy landscapes} in material sciences, partly to the re-discovery of the
{\em interpretational} difficulties of quantum mechanics in recent quantum optics.

\subsection{Theory -- model -- experiment} Physics offers an extremely useful practical
distinction between {\em theory}, {\em model} and {\em experiments}.
From his deep insight in astronomy, computing, linguistics and psychology,
Peter Naur ridicules such distinctions as ``metaphysical exaggeration" in \cite{Nau:KML}.
He may be right. We certainly should not exaggerate the distinction.
In this review, however, the distinction helps to focus on differences of the role of mathematics
in doing science.

\subsubsection{First principles} By definition, the very core of modelling is mathematics.
Moreover, if alone by the stochastic character of observations, but also due to the need to
understand the mathematics of all transducers involved in measurements, mathematics has its
firm stand with experiments. First principles, however, have a different status:
they do not earn their authority from the elegance of being mathematically wrapped, but from the almost infinite
repetition of similar and, as well disparate observations connected to the same
principle(s). In the first principles, mathematics and physics meet almost on eye level:
first principles are also {\em established} - like mathematics, and are only
marginally questioned. To me, the problem with the pretended eternal authority
of first principles is that new cosmological work indicates that the laws of nature
may also have undergone some development; that there might have ``survived" some
evolutionary relicts; and that we had better be prepared to be confronted
under extreme experimental conditions, with phenomena and relations which fall
out of the range of accredited first principles. The canonical candidate for such
a relict is the Higgs particle, whether already observed or not, see Holger
Bech Nielsen's contributions on the Quantum Gravity Assessment Workshop 2008, \cite{BooEspLes:QG}.

\subsubsection{Towards a taxonomy of models}Not necessarily for the credibility of
mathematical models, but for the way to check the range of credibility, the
following taxonomy of models may be extremely useful.

The Closing Round Table of the International Congress of
Mathematicians (Mad\-rid, August 22-29, 2006) was devoted to the
topic {\em Are pure and applied mathematics drifting apart?} As
panellist, Yuri Manin subdivided the mathematization, i.e., the way
mathematics can tell us something about the external world, into
three modes of functioning (similarly Bohle, Boo{\ss} and Jensen
1983, \cite{BoBoJe:IEO}, see also \cite{Bo:AIF}):



An {\em (ad-hoc, empirically based) mathematical model} ``describes
a certain range of phenomena, qualitatively or quantitatively, but
feels uneasy pretending to be something more." Manin gives two
examples for the predictive power of such models, Ptolemy's model of
epicycles describing planetary motions of about 150 BCE, and the
standard model of around 1960 describing the interaction of
elementary particles, besides legions of ad-hoc models which hide
lack of understanding behind a more or less elaborated mathematical
formalism of organizing available data.

A {\em mathematically formulated theory} is distinguished from an
ad-hoc model primarily by its ``higher aspirations. A theory, so to
speak, is an aristocratic model." Theoretically substantiated
models, such as Newton's mechanics, are not necessarily more precise
than ad-hoc models; the coding of experience in the form of a
theory, however, allows a more flexible use of the model, since its
embedding in a theory universe permits a theoretical check of at
least some of its assumptions. A theoretical assessment of the
precision and of possible deviations of the model can be based on
the underlying theory.

A {\em mathematical metaphor} postulates that ``some complex range
of phenomena might be compared to a mathematical construction". As
an example, Manin mentions artificial intelligence with its ``very
complex systems which are processing information because we have
constructed them, and we are trying to compare them with the human
brain, which we do not understand very well -- we do not understand
almost at all. So at the moment it is a very interesting
mathematical metaphor, and what it allows us to do mostly is to sort
of cut out our wrong assumptions. If we start comparing them with
some very well-known reality, it turns out that they would not
work."


Clearly, Manin noted the deceptive formal similarity of the three
ways of mathematization which are radically different with respect
to their empirical foundation and scientific status. He expressed
concern about the lack of distinction and how that may ``influence
our value systems". In the words of \cite[p. 73]{Bo:AIF}:

\citat{Well founded applied mathematics generates prestige which is
inappropriately generalized to support these quite different
applications. The clarity and precision of the mathematical
derivations here are in sharp contrast to the uncertainty of the
underlying relations assumed. In fact, similarity of the
mathematical formalism involved tends to mask the differences in the
scientific extra-mathematical status, in the credibility of the
conclusions and in appropriate ways of checking assumptions and
results... Mathematization can -- and therein lays its success --
make existing rationality transparent; mathematization cannot
introduce rationality to a system where it is absent... or
compensate for a deficit of knowledge.}

\noi Asked whether the last 30 years of mathematics' consolidation
raise the chance of consolidation also in phenomenologically and
metaphorically expanding sciences, Manin hesitated to use such
simplistic terms. He recalled the notion of Kolmogorov complexity of
a piece of information, which is, roughly speaking,

\citat{the length of the shortest programme, which can be then used
to generate this piece of information... Classical laws of physics
-- such phantastic laws as Newton's law of gravity and Einstein's
equations -- are extremely short programmes to generate a lot of
descriptions of real physical world situations. I am not at all sure
that Kolmogorov's complexity of data that were uncovered by, say,
genetics in the human genome project, or even modern cosmology data
... is sufficiently small that they can be really grasped by the
human mind.}


\section{Medicine}

\subsection{The special place of medicine}
Physicists of our time like to date the physics' beginning back to
Galileo Galilei and his translation of measurable times and
´distances on a skew plane into an abstract fall law. Before Galilei
- and long time after him, the methodological scientific status of
what we would call mechanical physics was quite low as compared with
medicine. Physics was a purely empirical subject. It was about
precise series of observations and quantitative extrapolations. It
was the way to predict planetary positions, in particular eclipse
times, the content of silver in compounds, or the manpower required
to lift a given weight with given weight arm. It was accompanied and
mixed up with all kinds of speculations about the spirits
and ghosts at work. We can easily see the continuity of results, of
observations and calculations from Kepler and Newton to our time.
However, we can hardly recognize anything in their thinking about
physics, in the way they connected physics with cosmic music or
alchemy.

\subsubsection{The maturity of medicine}
Contrary to that, from the rich ancient literature preserved, see
Diepgen \cite{Die:GME}, Kudlien \cite{Kud:BMD} and, in particular,
J{\"u}rss \cite[312--315]{Jur:GWD}, we can see that the mind set in
Greek medicine already from the fifth century BCE was {\em ours}:
instead of the partition (familiar from earlier and shaman medicine
and similar to the mind set preserved, as seen above, in physics
until recent times) into an empirical-rational branch (healing
wounds) and a religious-magic branch (cure inner diseases), a
physiological concept emerged which focused on the patient as an
individual organism within a population, with organs, liquids and
tissue, subjected to environmental and dietetic influences and, in
principle, open for unconfined investigation of functions, causal
relations and the processive course of diseases. In Hippocratic
medicine, we meet for the first time the visible endeavour after a
rational surmounting of all problems related to body events.

With a shake of the head, we may read of Greek emphasis and
speculations about the body's four liquids or other strange things,
like when we recall today the verdict of the medical profession 60
years ago against drinking water after doing sports and under
diarrhoea or their blind trust in antibiotics, not considering
resistance aspects at all. Contrary to physics, we have no
continuity of results in medicine, but, also contrary to physics, we
have an outspoken continuity in mind set: no ghosts, no metaphysical
spirits are permitted to enter our explanations, diagnoses,
prevention, curation and palliation.

To me, medicine is a science that is characterized by maturity in thinking -
but, until now, of ephemeral quality
of results. There are
astonishing progresses in the mathematization of biology and
medicine of the last decades. There are good reasons to expect that
medicine will become the most important field of mathematics
application - as physics has been in the last centuries.

\subsubsection{Hopeless aspects of medicine mathematization}
Before proceeding with this review, three warnings shall be presented.
The first warning is taken from Jean le Rond D'Alembert, Discours pr{\'e}liminaire
to the Encyclop{\'e}die, \cite[page vii]{dAl:DPE}:

\citat{It must be confessed, however, that geometers sometimes abuse
this application of algebra to physics. Lacking appropriate
experiments as a basis for their calculations, they permit
themselves to use hypotheses which are most convenient, to be sure,
but often very far removed from what really exists in Nature. Some
have tried to reduce even the art of curing to calculations, and the
human body, that most complicated machine, has been treated by our
algebraic doctors as if it were the simplest or the easiest one to
reduce to its component parts. It is a curious thing to see these
authors solve with the stroke of a pen problems of hydraulics and
statics capable of occupying the greatest geometers for a whole
lifetime. As for us who are wiser or more timid, let us be content
to view most of these calculations and vague suppositions as
intellectual games to which Nature is not obliged to conform, and
let us conclude that the single true method of philosophizing as
physical scientists consists either in the application of
mathematical analysis to experiments, or in observation alone,
enlightened by the spirit of method, aided sometimes by conjectures
when they can furnish some insights, but rigidly dissociated from
any arbitrary hypotheses.}

\noi This pompous verdict has not kept d'Alembert from contributing nine years later a mathematical
model to smallpox inoculation epidemiology which even today is worth reading, see
Dietz and Heesterbeek \cite[Section 7]{DieHee:DBE}.

A second warning against exaggerated expectations regarding the mathematization
of biology and medicine comes from the suspicion that many of the organs, organelles
and DNA sequence regions are {\em evolutionary relicts} without any active function and
with the only meaning to confuse the observer and the modeller.

A third warning relates to the conflict between {\em reductionist} and {\em holistic}
approaches. We have seen a similar conflict in physics when discussing the chances vs. vanity 
of the unifying approaches dating back to Einstein and re-actualized in various approaches to quantum gravity.
For medicine it seems clear that a strictly reductionist program is mandatory when we wish to replace
ad-hoc assumptions by reference to first principles. On the other side, there is no doubt that
most body functions and processes involve many parts of a cell, organ, and organism at the same time.
Understandably, the slogan of holistic {\em Systems Biology} has become very popular and great
expectations are attached to it.
Both programs will disclose new and interesting facts and features and
offer mathematicians rich fields to work in. To me, however,
the most promising approach will be somewhere in between. Perhaps {\em focussed} systems biology
will prove to be able to hit the wall and make a hole in the wall,
i.e., a breakthrough. That is a medicine and biology discarding many surely relevant aspects, lengths, distances, and focussing
on a restricted, but multi-level and multi-scale set of processes by applying a wide range of mathematical methods and first physics principles. For an example see below Paragraph
\ref{ss:cells}.

\subsection{Populations, organisms, organs, cells, DNA}
In this report, we shall try to separate these levels wherever possible. On each level,
quite different types of mathematical methods are applied and quite different types of results obtained.

\subsubsection{Populations} 
The intersection between humans and infectious diseases is perhaps the most intensively
studied problem in population biology (partly originating in the actuarial literature,
see the informative study \cite{DieHee:DBE}).
Epidemiological studies about the distribution of disease in human populations
and the factors determining that distribution, have shown to be useful
in public health and preventive medicine, in particular
when not only the interference was based on statistical methods but also the design of the
data collection and the control of the reliability of the reported data.

Perhaps the most important result of mathematical epidemiology has been the statistical proof of
the health risks of smoking - in spite of the great statistician Ronald A. Fisher's insistence
that excess health problems and mortality of smokers might have two other clever explanations, besides the common,
namely that subconscious feeling of early death might be the main reason of the desire to smoke or
that a predisposition to lung cancer and arterial diseases are genetically confounded with the desire to
smoke, see \cite{Fis:CCA} and recognized and ridiculed the evidence as we would with the observed
correlation between the occurrence of storks and the frequency of birth in the Mark Brandenburg
in the 1920s (German children are told that the stork brings the babies in the same way as Anglo-Saxonic
babies are found under gooseberry bushes). It must be acknowledged that the work of R.A. Fisher and his colleagues
on correlation and regression analysis, goodness of fit, sampling and statistical error
in the first half of the 20th century had given epidemiology new techniques of measurement
and also undermined erroneous certainties. However, it remains strange that Fisher
in his polemic did not disclose his role as advisor to the Tobacco Manufacturers Standing
Committee nor admitted that a less libertarian world view would support the statistical results
and health recommendations by Sir Richard Doll, the discoverer of the smoking correlation.

While main stream epidemiology has its greatest successes in
identifying environmental and behavioural factors, {\em genetic
epidemiology} impresses by identifying hereditary factors, in
particular providing large panels of the most relevant single
nucleotide polymorphisms (SNPs) associating in selected diseases.
Recently, a SNP whole genome scan has been carried out on 1599
patients with type 2 diabetes and 1503 control subjects matched from
Finland and Sweden, mapping loci for various elements of metabolic
syndrome, see also Groop et al. \cite{Gro:GWA}. While genotyping now
begins to become rather straight forward, it often remains unclear
how to translate the new information into theoretical analysis of
the disease and to design new strategies for early diagnosis
(perhaps not at all always desirable) and treatment, see also below
Paragraph \ref{ss:cells}.

Another branch of epidemiology are the SIR models for {\em
infectious diseases}. For teaching qualitative reasoning and
numerical analysis of coupled non-linear differential equation, the
modelling of {\underline S}usceptibles, {\underline I}nfectives and
{\underline R}emoved, respectively, {\underline R}esistents have
been a beloved dull for decades: a toy of pure mathematics and
applied terminology like a sheep in wolf's clothing.

The power of SIR models became evident, however, in controversies 
about vaccination schemes against Rubella, where the models
could predict that non-vaccination would be more safe than
non-compulsory vaccination below the critical level of coverage and
Greece had to pay for not listening to the predictions, \cite[pp.
77-79, 102-110, 152f, 167f, 188-199]{AndMay:IDH}. Modified SIR
models have also proved applicable for estimating various factors
for the spread (respectively, control) of Methicillin Resistant {\em
Staphylococcus aurea} (MRSA): the role of antibiotic consumption,
the role of isolation and other sanitary procedures in hospitals,
the comparative weights of entrance and exit screening to keep the
community pool down, etc., see \cite{BBD:CMR}.

It seems that long-time oscillations and rapidly mutating vira like the HIV virus
or the
bird flu H5N1 still challenge mathematicians and evolutionary epidemiologists 
due to the non-linar aspects of the processes. Moreover, it seems necessary
and possible by advances in computing, to replace the SIR models by time-series
SIR (TSIR). Here the point is to drop the homogeneous-mixing assumption and to
admit that not all hosts have identical rates of disease-causing contacts, see
Grenfell \& collaborators \cite{BanGreMey:WIB}.

\subsubsection{Organisms and organs}
To give one (outstanding) example, we refer to the modelling of the
cardio-vascular system, \cite{OttDan:MMM}. As discussed above in
paragraph \ref{ss:prediction}, computer supported modelling can
provide the most beautiful animations. However, modelling can have
other goals: explaining strange, but well-established phenomena
(e.g., that the forces exercised by the heart are slightly below the
forces measured by the blood stream); supporting heart surgery and
implantations; and tuning anaesthesia simulators for the
training of anaesthesiologists. For these goals, the reliability of
the modelling's basic assumptions has shown to be decisive. Compared
with the most smashing animations, modelling the cardio-vascular
system like a sewage system with elastic walls may sound primitive.
However, it gains its reliability by its firm foundation on first
principles.

\subsubsection{Cells}\label{ss:cells}
The mathematization of the cell has many levels and many scales. To
give an example, I shall describe an evolving - focussed - systems
biology of regulated exocytosis in pancreatic $\beta$-cells, mostly
based on \cite{Boo:REX}. These cells are responsible for the
appropriate insulin secretion. Insufficient mass or
function of these cells characterize Type 1 and Type 2 {\em diabetes
mellitus} (T1DM, T2DM). Similar secretion processes happen in nerve
cells. However, characteristic times for insulin secretion are
between 5 and 30 minutes, while the secretion of neurotransmitters
is in the millisecond range. Moreover, the length of a $\beta$-cell
is hardly exceeding 3000 nanometres (nm), while nerve cells have
characteristic lengths in the cm and meter range. So, processes in
$\beta$-cells are more easy to observe than processes in nerve
cells, but they are basically comparable.

It seems that comprehensive research on $\beta$-cell function and mass has been seriously hampered for 80
years because of the high efficiency of the symptomatic treatment of T1DM and T2DM by insulin injection.
Recent advances - and promises - of nanotechnology suggest the following radically new research agenda, to be
executed first on cell lines, then on cell tissue of selected rodents, finally on living human cells:



\noi 1. Synthesize magneto-luminescent nanoparticles; develop a
precisely working electric device, which is able to generate a
properly behaving electromagnetic field; measure cytoskeletal
viscosity and detect the interaction with organelles and actin
filaments by optical tracking of the forced movement of the
nanoparticles. Difficulties to overcome: protect against protein
adsorption by suitable coating of the particles and determine the
field strength necessary to distinguish the forced movement from the
underlying Brownian motion.


\noi 2. Synthesize luminescent nanoparticles with after-glow
property (extended duration of luminescence and separation of
excitation and light emission); dope the nanoparticles with suitable
antigens and attach them to selected organelles to track
intra-cellular dynamics of the insulin granules.


\noi 3. Develop a multipurpose sensor chip and measure all electric
phenomena (varying potential over the plasma membrane, the bursts of
$Ca^{++}$ ion oscillations,  and changing impedances on the surface
of the plasma membrane for precise chronical order of relevant
secretion events.

\noi 4. Describe the details of the bilayer membrane-granule fusion
event (with the hard numerical problems of the meso scale, largely
exceeding the well-functioning scales of molecular dynamics).


\noi 5. Connect the preceding dynamic and geometric data with
reaction-diffusion data.

Connect the preceding dynamic and geometric data with genetic data.

Develop clinical and pharmaceutical applications:
\begin{itemize}
\item Quality control of transplants for T1DM patients.
\item Testing drug components for $\beta$-cell repair.
\item Testing nanotoxicity and drug components for various cell types.
\item Early in-vivo diagnosis by enhanced gastroscopy.
\item Develop mild forms of gene therapy for patients with  over-expressed  major type 2-diabetes gene TCF7L2
by targeting short interfering RNA sequences (siRNAs) to the
$\beta$-cells, leading to degradation of excess mRNA transcript.
(This strategy may be difficult to implement, due to the degradation
of free RNA in the blood and the risk of off-target effects.)
\end{itemize}

We shall not go into the mathematical details of the involved compartment models, free boundary theory,
reaction-diffusion equations,
data analysis etc. Summing up, the mathematization can help design
relevant experiments and support the determination of basic parameters.
Mathematical modelling and simulation can point to {\em
possible} new phenomena or new relations which have to be confirmed and investigated more thoroughly in the
laboratory. However, at the present state of our knowledge about cell and cell membrane processes,
one must doubt whether making predictions {\em in silico} can really
replace experiments {\em in vivo} or {\em in vitro}.


\section{Economics}

While a holistic unifying view in physics has a smell of vanity (in spite of efforts in quantum
gravity) and the time has hardly come for a holistic all-embracing systems biology programme
in medicine, a rational point of departure for economics, in particular under the present crisis,
can only be a {\em systems view}.

On the one side we have the relief of breaking the growth curve which
has been depriving future generations of non-renewable resources. In
some parts of the world we can say ``Enough is enough". We can choose
a conscious way of living, we can replace purely quantitative scales
for the quality of life, like GNP, money, work places, hospital beds
and spending power by other social and individual values like
health, peace, prosperity, happiness. We can look to Cuba or other
model countries who have shown how to live with dramatically
restricted energy consumption - if a lot more should go wrong.

On the other side we must worry that the crisis, like the Great
Depression of 1929-1945, see \cite{Kin:TWD} will deepen the
differences between rich and poor within and between countries, and
can foster unprecedented political and military catastrophes.

Many scientific subjects are challenged to give explanation and
advice in the present crisis. We would be grateful for relief in
understanding what is going on,  and qualified warning would be
needed against leaving the required decisions to the isolated world
of so-called {\em decision makers} or {\em chief advisors} and {\em
chief economists}. Perhaps not mathematical economics is needed
(see, however, the recent ex-post analyses \cite{CCFS:FAD}),
but rather psychology which is proficient in how individuals and
masses experience crises, or philosophy which should be proficient
on the relation between property and lack of sense of
responsibility.

What can a mathematician add to the avalanche of crisis comments in the literature?
Where are the pitfalls and sins and where the chances of mathematical modelling in economics and
finance under the present crisis?

\subsection{Praise and criticism of greed and accountability}
In laborious analysis of hundred thousands of figures and tables,
the classic socialist thinkers, culminating with Friedrich Engels
and Karl Marx, have shown that capitalism builds on, supports and
develops {\em Rechenhaftigkeit -- accountability}, and that abstract
greed, greed for more money, greed for more profit, is the most
powerful and progressive aspect of capitalism. They advocated for
another society where decisions, human relations, development should
not be driven by greed but by conscious and debated choices. That
makes sense. But clinging to capitalism, though preferably without
greed, is media hype and logical non-sense.

One does not need relatively advanced macro-economical predictive
models to see that greed, while driving economic growth and
economisation of (many, not all) resources, generates
contradictions. It generates abundance of free capital in the
presence of unsatisfied consumption wishes and deepens the social
divide. Under this perspective, it was not at all unreasonable by
the US-American Federal Bank to admit the application of {\em
subprime} loans of questionable security: many young men of
US-American lower middle class were risking health and life in Iraq
and Afghanistan. As Alan Greenspan confessed: under the given
circumstances, the subprime loans were necessary. It {\em gave more
people a stake in the future of our country and boded well for the
cohesion of the nation} in a country at war, \cite[Chapter
11]{Gre:AT}!

In many comments, the so-called {\em neo-classical economic equilibrium theory} is blamed for
having closed the eyes of politicians and bankers in front of the upcoming crisis. Mathematical
economics can, in deed, be blamed for many misleading concepts, statements and predictions. Mathematical
models have had their share in moving the focus from basic contradictions in the sectors of production
and consumption (private daily and luxury requirements and state and military expenditures,
i.e., {\em Main Street}) over into the sphere  of circulation ({\em Wall Street}) proliferating
for much too long the impression that the financial crisis was only a problem for banks and pension funds.
No doubt, floods of money and signs are more easy to quantify and mathematize than the state and
development of productive forces and production relations.

However, mathematical models can not really be blamed for the crisis. The crisis was
there, it was open, it was visible, and one reacted correctly according to the conditions of
capitalism - and made so the crisis to deepen and to propagate into a broader crisis.

Now, the system-immanent rational continues when greed and state power are merged in an unprecedented way
and huge and almost free self-service tables are provided for poisoned capital on the verge of bankruptcy.
Once again, it's clearly the correct answer within a context which is considered as fixed. Consequently,
one must be afraid that these help packages will further deepen the crisis, at least for a time to come
until the financial markets will regain a certain transparency and will show their force to clean-up the mess
- with sacrifices of yet unknown magnitude and composition to be offered to the gods.

\subsection{Hedging and the fight for security}
Mathematical modelling in finance and economics has not only greed, growth and profit as
variables to predict and optimize, but also the security of investments, the safety of an
portfolio of private investors, hedge funds and pension funds. Actually, mathematicians, actuaries
and mathematically trained economists can make calculations that appear valid, credible and trustworthy.
Thirty years ago, there were only two mathematicians employed by finance in the whole of Denmark. Since
a couple of years, often more than half of a year of graduates in math and physics end up there.

Credit evaluation and fixing prices of new financial creations is so intricate that parts of
the most modern mathematics, e.g., martingales and heat equation, are challenged. In {\em normal}
situations, the calculations hold. That's the way they are {\em calibrated}. These calculations have had,
admittedly, their part in conjuring false security. It is ironical that the mathematicians were employed
for enhancing the security of investments by modelling. However, by making the trade with
options and other financial derivates more just, more reliable, more easy, we have also created a machine
that can magnify and worsen a developing crisis by market automatics.
See also Krickeberg \cite{Kri:MEM} on the one hand, and on the other Carmona \& Sircar \cite{CarSir:MFC}.

\subsection{Soft selection instead of control of the irrational?}
In discussions about the crisis one often hears the words {\em
complexity}, {\em uncertainty}, {\em control} and {\em chaos}. In
mathematics we know many very simple, even deterministic systems
that are absolutely unpredictable as long as they are permitted to
follow their own way. Take the double pendulum: it can turn
somersaults. Small, practically invisible differences decide whether
there comes one somersault more. Unpredictable? Yes. Uncontrollable?
Yes - unless we put a sufficiently narrow box around the pendulum.
Then no further somersaults will come. Perhaps such a box should be
built for our economic systems as well, sooner or later.

As discussed in the previous chapters, mathematical models are,
perhaps like all science, able to predict concrete events (like eclipses
and election results) and consequences of performed operations
(e.g., when we put on the light at the switch). Mathematics and
science, however, can do that in ideally-simple cases only. In
complex cases, we will be lucky if we at least can point to
necessary changes, confinements, regulations to create a predictable
and adjustable system. Instead of the possibly vain goals of
controlling an irrational system we may consider to soften our
system of goal functions. In the long run, soft selection and
redesign of priorities and ways of cooperation, as explained
in \cite{Gal:ESS}, will be the most efficient way to peace and
prosperity. To me, that is the only proved theorem of this paper.

\section{General trends of mathematization and modelling}

\subsection{Deep divide}\label{ss:dd}
Regarding the power and the value of mathematization, there is a deep moral divide
both within the mathematics community and the public.

On the one side, we have the outspoken science and math optimism of
outstanding thinkers: Henri Poincar{\'e}'s {\em Nature not only
suggests to us problems, she suggests their solution}; David
Hilbert's {\em Wir m{\"u}ssen wissen; wir werden wissen - We must
know; we will know} of his Speech in K{\"o}nigsberg in 1930, now on
his tomb in G{\"o}ttingen; or Bertolt Brecht's vision of
mathematical accountability in {\em Die Tage der Kommune}
\cite{Bre:DTK} of 1945: ``Das ist die Kommune, das ist die
Wissenschaft, das neue Jahrtausend... - That is the Commune, that is
the science, the new millennium..."). We have astonishing evidence
that many mathematization concepts either appear to us as {\em
natural} and {\em a-priori}, or they use to emerge as clear over
time. We have the power and validity of extremely simple concepts,
as in {\em dimension analysis}, {\em consistency requirements} and
{\em gauge invariance} of mathematical physics. Progressive
movements emphasize science and education in liberation movements
and developing countries. Humanitarian organizations (like WHO and
UNICEF) preach science and technology optimism in confronting mass
poverty and epidemics.

On the other side, deep limitation layers of science and
mathematical thinking have been dogged up by Kurt G{\"o}del's {\em
Incompleteness Theorem} for sufficiently rich arithmetic systems,
Andrei N. Kolmogorov's {\em Complexity Theory}, and Niels Bohr's
notion of {\em Complementarity}. Incomprehensibility and lack of
regularity continue to hamper trustworthy mathematization. Peter Lax
\cite[p. 142]{Lax:MAP} writes about the {\em profound mystery of
fluids}, though recognizing that different approaches lead to
remarkable coincidence results, supporting reliability.

The abstruseness of the mathematical triumphs of the hydrogen bomb
is commonplace. The wide-spread trust in superiority and
invincibility, based on mathematical war technology like high
precision bombing, has proved to be even more vicious for warriors
and victims than the immediate physical impact of the very
math-based weaponry, recently also in Iraq and Afghanistan.

In between the two extremes, we have the optimistic scepticism of
Eugene Wigner's {\em unreasonable effectiveness of mathematics}, but
also Jacob Schwartz's verdict against the {\em pernicious influence of mathematics on science}
and Albert Einstein's demand
for {\em finding the central questions} against the dominance of the beautiful and the difficult.



\subsection{Charles Sanders Peirce's semiotic view}
From the times of Niels Bohr, many physicists, mathematicians and
biologists have been attentive to philosophical aspects of our
doing. Most of us are convinced that the frontier situation of our
research can point to aspects of some philosophical relevance - if
only the professional philosophers would take the necessary time to
become familiar with our thinking. Seldom, however, we read
something of the philosophers which can inspire us.

The US-American philosopher Charles Sanders Peirce (1839-1914) is an admirable exception.
In his semiotics and pragmaticist (he avoids the word ``pragmatic") thinking, he provides
a wealth of ideas, spread over an immense life work.
It seems to me that many of his ideas, comments, and concepts can shed light on the why and how of
mathematization. Here I shall only refer some thoughts of Peirce's {\em The Fixation of Belief}
from 1877, see \cite{Pei:FB}.

My fascination of Peirce's text is, in particular, based on the
following observations which may appear trivial (or known from
Friedrich Engels), but are necessary to repeat many times for the
new-modeller:


\noi 1. For good and bad, we are all equipped with innate (or
spontaneous) orientation, sometimes to exploit, sometimes to subdue.
Our innate orientation is similar to the  habits of animals in our
familiar neighbourhood. We are all "logical machines".


\noi 2. However, inborn logic is not sufficient in foreign (new)
situations. For such situations, we need methods how to fixate our
beliefs. Peirce distinguishes four different methods. All four have
mathematical aspects and are common in mathematical modelling.
\begin{description}
\item[Tenacity] is our strength not to become confused, not to be blown away by unfounded arguments, superficial objections, misleading examples, though sometimes keeping our ears locked for too long.

\item[Authority] of well-established theories and results is what we tend to believe in
and have to stick to. We will seldom drop a mastered approach in favour of something new and unproved.

\item[Discussion] can hardly help to overcome a belief built on tenacity or authority.

\item[Consequences] have to be investigated in all modelling. At the end of the day,
they decide whether we become convinced of the validity of our approach (Peirce's Pragmaticist Maxim).

\end{description}


\noi 3. The main tool of modelling (i.e., the fixation of belief by
mathematical arguments) is the transformation of symbols (signals,
observations, {\em segments of reality}) into a new set of symbols
(mathematical equations, {\em models} and {\em descriptions}). The
advantage for the modeller, for the person to interpret the signs,
is that signs which are hard or humid and difficult to collect in
one hand can be replaced by signs which we can write and manipulate.


\noi 4. The common mapping cycle {\em reality} $\to$ {\em model}
$\to$ {\em validation} is misleading. The quality of a mathematical
model is not how similar it is to the segment of reality under
consideration, but whether it provides a flexible and goal-oriented
approach, opening for doubts and indicating ways for the removal of
doubts (later trivialized by Popper's falsification claim). More
precisely, Peirce claims


Be aware of differences between different approaches!


Try to distinguish different goals (different priorities) of
modelling as precise as possible!


Investigate whether different goals are mutually compatible, i.e.,
can be reached simultaneously!


Behave realistically! Don't ask: {\em How well does the model reflect a
given segment of the world}? But ask: {\em Does this model of a
given segment of the world support the wanted and possibly wider
activities / goals better than other models}?








\subsection{Integral geometry, for ever}
Perhaps a brief summary of the history of {\em integral geometry},
from Cavalieri's Principle to Harmonic Analysis, Characteristic Classes, Point Processes, and Stochastic Geometry
could illustrate some of Peirce's points and explain the origin and development of far-reaching
mathematical concepts, the power of mathematization, from gauge theoretic physics to the counting of
brain cells, and the continuing and deepening limitations of mathematization over urging problems,
be it the limitations of
mathematical concepts like {\em homotopy invariance} (when confronted with spectral invariants)
or limitations due to peculiarities of the subject domain. There is more to discuss.



\end{document}